# An Alternative method in Multi-Attribute Decision Making using Data Envelopment Analysis and Fuzzy concept

Majid Zerafat Angiz L.[*1], Mohd Kamal Nawawi[1], Mohammad Ghadiri[2], Adli Mustafa[3]

[*1]School of Quantitative Sciences, Universiti Utara Malaysia, 06010 Sintok, Kedah, Malaysia

[2]Networking & Computer Company, 1820 San Pedro Dr NE Albuquerque, New Mexico, 87110, United States

[3]School of Mathematical Sciences, Universiti Sains Malaysia, Penang, Malaysia

* Corresponding author.  Tel.: +17637770788

E-mail address: mzerafat24@yahoo.com (M. Zerafat Angiz L.)

**Abstract**

Data Envelopment Analysis (DEA) as mathematical models evaluates the technical efficiency of Decision Making Units (DMU) having multiple inputs and multiple outputs. Researchers are interested in applying DEA models in Multi Attribute Decision Making (MADM) environment, but evaluation by these models is different in nature than MADM. This is why the results are not satisfactory. In this paper first, a challenging discussion is provided to indicate ranking using traditional DEA models is not reliable, and then a hybrid model using DEA and fuzzy concepts is proposed to present a self-assessment for each DMU.

## 1. Introduction

Data envelopment analysis (DEA) is a productivity analysis tool that estimates production frontiers and measures the relative efficiency of decision making units (DMUs). In other words, DEA measures technical efficiency of DMUs. In this view, in input orientation, the DMU under evaluation is efficient if it produces the specific outputs with the consumed specific inputs, and simultaneously, there is no other DMU consuming the less value of inputs, compared with the DMU under evaluation, and produce the same value of outputs. In output orientation, the DMU under evaluation is efficient if it produces the specific value of outputs with the consumed specific value of inputs, and simultaneously, there is no other

DMU to consume the same value of inputs compared with the DMU under evaluation and produce the more value of outputs. After publishing the first DEA paper proposed by Charnes et al. [1], known as the CCR model, Banker et al. [2] developed a variable returns to scale variation of the model called BCC model. Both models have undergone various theoretical extensions and also many successful applications [3]. The CCR and BCC models appraise the radial efficiency but do not consider non-zero slacks. This is why The CCR model evaluates the DMUs on the weak efficient frontier as efficient. To make up to this deficiency, the additive model presented in [4], based on non-radial measure, was introduced whose objective function is summation of the slacks. This model can discriminate the efficient and inefficient DMUs. The *CCR*, *BCC* and ADD models are the main models in DEA, but they do not exhaust the available DEA models. Another non-radial approach which deals with slacks directly is a slack-based measure of efficiency (SBM) introduced by Tone [5].

According to Adler et al. [3], the existed ranking methods have been divided into six categories; the cross-efficiency methods initiated by Sexton et al.[6], super efficiency based approaches pioneered by Andersen and Peterson [7], benchmarking ranking methods, the ranking method with multivariate statistics in the DEA models and the methods based on ingredient of DEA and multi-criteria decision making methods.

An alternative technique based on non-radial measure is applied to discriminate efficient DMUs in [8] and subsequently Saati et al. [9] modified the non-radial model presented in [8] and converted it to an input-output orientation model that caused the LP model to be always feasible. In this method, the traditional DEA models are reformulated by excluding the input and output of the DMU that is being ranked from the models

There are several hybrid models in the DEA literature. In some models, other methodologies are utilized to support DEA models to obtain more appropriate results. Sinuany-Stern *et al.* [10] proposed an approach based on the relationship between DEA and AHP (analytic hierarchy process) to rank DMUs. In another hybrid measure approach, Tone [5] introduced a model in which the radial and non-radial measure approaches were simultaneously integrated into the DEA mathematical program.

The performance evaluation models presented in the literature on DEA are divided to two categories. In the first category, the methods retouch the production possibility set (PPS) to provide special aims. The super efficiency models and models based on return to scale are in this category. In the second, other models aim DMUs in the same efficiency frontier using different measures. For instance, the efficient

DMU in the efficiency frontier corresponding to an inefficient DMU under evaluation in the CCR model and additive model are different. Therefore, two different optimal solutions are obtained.

In contrary to previous methods, in this paper an alternative approach is presented in which production frontier of each DMU is different from others. We concentrate on the opposite side of production possibility set corresponding with the DMU under evaluation. In other words, we assume that there is just one DMU and so, the production frontier passes through the DMU under evaluation. Based on this assumption, a hybrid model basis on the super efficiency model presented in [9], and fuzzy interpretation of production possibility set, is introduced. Notice that the inputs and outputs are crisp data.

The remainder of this paper is organized as follows: In Section 2, some background information about fuzzy numbers and the CCR model are given. A new approach is described in Section 3. In Section 4, a case study is given. Finally, the conclusion is presented in Section 5.

**2. Background**

**2.1 Data envelopment analysis**

In order to review the DEA models, the following two general assumptions are specified:

**a)** There are $n$ *DMU*s denoted by $j \in J$, each of which produces a nonzero output vector $Y_j = (y_{1j}, y_{2j}, ..., y_{sj})^t \geq 0$ using a nonzero input vector $X_j = (x_{1j}, x_{2j}, ..., x_{mj})^t \geq 0$, where the superscript $'t'$ indicates the transpose of a vector. Here, the symbol $'\geq'$ indicates that at least one component of $X_j$ or $Y_j$ is positive while the remaining $X_j$'s or $Y_j$'s are nonnegative.

**b)** There is no *DMU* in $J$ whose data domain can be proportionally expressed by that of another *DMU*.

**Definition 2.1.** Given the (empirical) points $(X_j, Y_j)$, $j = 1, 2, ..., n$, the Production Possibility Set (PPS) is defined as follows:

$$T = \{(X_t, Y_t) \mid \text{output } Y_t \text{ can be produced by input } X_t\}$$

**Definition 2.2.**

**a)** The production possibility $(X_t, Y_t)$ is a frontier point (input-oriented) if $(\alpha X_t, Y_t) \in T$ implies $\alpha \geq 1$.

**b)** Production possibility $(X_t, Y_t)$ is a frontier point (output-oriented) if $(X_t, \beta Y_t) \in T$ implies $\beta \leq 1$.

**a)** The production possibility $(X_t, Y_t)$ is a frontier point (input-oriented) if $(\alpha X_t, Y_t) \in T$ implies $\alpha \geq 1$.

**b)** Production possibility $(X_t, Y_t)$ is a frontier point (output-oriented) if $(X_t, \beta Y_t) \in T$ implies $\beta \leq 1$.

To construct the production possibility set, the following postulates are assumed:

**1)** (Ray Unboundedness) If $(X_t, Y_t) \in T$ then $(\gamma X_t, \gamma Y_t) \in T$ for $\gamma > 0$.

**2)** (Convexity) If $(X_t, Y_t) \in T$ and $(X_u, Y_u) \in T$ then

$(\lambda X_t + (1-\lambda) X_u, \lambda Y_t + (1-\lambda) Y_u) \in T$ for all $\lambda \in [0,1]$.

**3)** (Monotonicity) If $(X_t, Y_t) \in T$, $X_u \geq X_t$ and $Y_u \leq Y_t$ then $(X_u, Y_u) \in T$.

**4)** (Inclusion of Observation) All the observations belong to production possibility set.

**5)** (Minimum Extrapolation) If $T'$ be a set different from $T$ which satisfies the mentioned above postulates, then $T \subseteq T'$.

The production possibility set corresponding to constant return to scale constructed with the aforementioned postulates will be as follows:

$$T_c = \left\{ (X_t, Y_t) / X_t = \sum_{j=1}^{n} \lambda_j x_{ij}, Y_t = \sum_{j=1}^{n} \lambda_j y_{rj}, \lambda_j \geq 0, \begin{array}{l} j=1,2,...,n \\ i=1,2,...,m \\ r=1,2,...,s \end{array} \right\} \quad (2.1)$$

Constant return to scale (CRS) means that an increase in the amount of inputs consumed leads to a proportional increase in the amount of outputs produced and if this increase is culminated in larger or smaller than proportional increase in the amount of outputs, return to scale will be increasing (IRS) or decreasing (DRS), respectively.

To evaluate efficiency corresponding to set $T_c$, consider the following model.

$$\theta_p^* = \min \theta_p \quad (2.2)$$
$$\text{s.t.}$$
$$(\theta_p X_p, Y_p) \in T$$

CCR model (input-oriented) for evaluating the efficiency of $DMU_p$, is written as follows:

$$\theta_p^* = \min \theta_p \quad (2.3)$$
$$\text{s.t.}$$
$$\sum_{j=1}^{n} \lambda_j x_{ij} \leq \theta_p x_{ip} \quad i = 1, 2, \ldots, m$$
$$\sum_{j=1}^{n} \lambda_j y_{rj} \geq y_{rp} \quad r = 1, 2, \ldots, s$$
$$\lambda_j \geq 0 \quad j = 1, 2, \ldots, n$$
$$\theta_p \text{ free}$$

The above CCR model is called as radial efficiency measures, because it optimizes all inputs or outputs of a DMU at a certain proportion. Färe and Lovell [5] introduced a non-radial measure which allows non-proportional reductions in positive inputs or augmentations in positive outputs. Saati et al. [9] suggested a non-radial model as follows:

$$\min \; \Omega_p = \omega_p + 1 \quad (2.4)$$
$$\text{s.t.}$$
$$\sum_{j=1}^{n} \lambda_j x'_{ij} \leq x'_{ip} + \omega_p 1 \quad \forall i$$
$$\sum_{j=1}^{n} \lambda_j y'_{rj} \geq y'_{rp} - \omega_p 1 \quad \forall r$$
$$\lambda_j \geq 0 \quad \forall j$$
$$\omega_p \text{ free}$$

In the above model, $x'_{ij} = \dfrac{x_{ij}}{\max\{x_{ij}\}}$, $y'_{rj} = \dfrac{y_{rj}}{\max\{y_{rj}\}}$ $\begin{cases} i=1,2,...m \\ r=1,2,...s \\ j=1,2,...n \end{cases}$, and $\omega_p$ is a free variable and measures the efficiency of $DMU_p$.

This model projects the DMU under evaluation on the frontier by decreasing the inputs and increasing the outputs.

## 2.2 Fuzzy number and Fuzzy Data Envelopment Analysis

### 2.2.1 Fuzzy number

A fuzzy number is a special case of a convex fuzzy set. In this section the L-R fuzzy number is explained.

**Definition1**. A fuzzy number $\tilde{x}$ is a convex normalized fuzzy set $\tilde{x}$ of the real line R such that

1. it exists exactly one $x_o \in R$ with $\mu_{\tilde{x}}(x_o) = 1$ ($x_o$ is called the mean value of $\tilde{x}$).

2. $\mu_{\tilde{x}}(x)$ is piecewise continuous.

Dubois and Prade [11] suggest a special case type of representation for fuzzy numbers of the following type:

**Definition2**. A fuzzy number $\tilde{x}$ is of LR-type if there exist reference functions L(for left), R (for right, and $\alpha > 0$, $\beta > 0$ with

$$\mu_{\tilde{x}}(\bar{x}) = \begin{cases} L\left(\dfrac{\bar{x} - x^l}{x^m - x^l}\right) & \text{for } x^l \leq \bar{x} \leq x^m \\ R\left(\dfrac{x^u - \bar{x}}{x^u - x^m}\right) & \text{for } x^m \leq \bar{x} \leq x^u \end{cases}$$

$x^m$, called the mean value of $\tilde{x}$, is a real number, and $\alpha = x^m - x^l$ and $\beta = x^u - x^m$ are called the left and right spreads, respectively, $x^l$ and $x^u$ are the lower value and the upper value of the interval of fuzzy number (see Figure 1). Symbolically, $\tilde{x}$ is denoted by $(x^m, x^l, x^u)$ or $(x^m, \alpha, \beta)$. If $\tilde{x}$ is a symmetric triangular fuzzy number, $R(x) = L(x) = x$ is implied.

### 2.2.2 Fuzzy CCR model

Due to lack of complete knowledge and information, precise mathematics is not sufficient to model a complex system. Although, in real world, decisions are based on qualitative as well as quantitative data, a fuzzy approach seems fit to deal with such problems.

CCR model has its production frontier spanned by the linear combination of the existing DMUs. But, production frontier in CCR and fuzzy CCR model are different. The frontiers of the CCR model have no flexibility and have linear characteristics, while those of the fuzzy are flexible.

$$\min \quad \phi_p \tag{2.5}$$
$$\text{s.t.}$$
$$\sum_{j=1}^{n} \lambda_j \tilde{x}_{ij} \leq \phi_p \tilde{x}_{ip}$$
$$\sum_{j=1}^{n} \lambda_j \tilde{y}_{rj} \geq \tilde{y}_{rp}$$
$$\lambda_j \geq 0$$
$$\phi_p \text{ free}$$

where "~" indicates the fuzziness.

As suggested in [8], a non radial CCR-DEA model with fuzzy coefficients is given as follows:

$$\min \quad W_p = w_p + 1$$
$$\text{s.t.}$$
$$\sum_{j=1}^{n} \lambda_j \tilde{x}_{ij} \leq \tilde{x}_{ip} + w_p 1 \quad \forall i \tag{2.6}$$
$$\sum_{j=1}^{n} \lambda_j \tilde{y}_{rj} \geq \tilde{y}_{rp} - w_p 1 \quad \forall r$$
$$\lambda_j \geq 0 \quad \forall j$$
$$w_p \text{ free}$$

### 2.2.3 The MADM methods

The two main categories for evaluating alternatives are available in the literature; compensatory and non-compensatory. In the Non-compensatory methods trade-off is not allowed. Dominance, Max-min, Max-max, Disjunctive-Satisfying and Lexicograph are categorized in this area. Among the above methods, Max-min, Max-max do not need information from decision makers. One of the most well-known method among compensatory category is TOPSIS which is compared with proposed model. Consider $n$ alternatives $A_j$ $(j = 1, 2, ..., n)$ and $m$ attributes $r_{ij}$ $(i = 1, 2, ..., m; r = 1, 2, ..., n)$.

### 2.2.3.1 Max-min

The Max-min means the maximum amount of profit. In this method, the scale of all attributes is converted in dimensionless form. For doing so, the following calculations are done;

$$h_{ij} = \frac{r_{ij}}{r_j^{max}} \quad i = 1, 2, ..., m \quad \text{for benifit attributes}$$

$$h_{ij} = \frac{r_j^{min}}{r_{ij}} \quad i = 1, 2, ..., m \quad \text{for cost attributes}$$

Where $r_j^{max} = \max_i \{r_{ij}\}$ $j = 1, 2, ..., n$ and $r_j^{min} = \min_i \{r_{ij}\}$ $j = 1, 2, ..., n$.

Most appropriate alternative is gained as follows:

$$A^* = \left\{ A_j \left| \max_j \min_i h_{ij} \right. \right\}$$

In order to use the Max-min methodology as a full ranking method, after selecting the best alternative, we eliminate it and redo the process of the normalization and search the next best one.

### 2.2.3.2 Characteristics of the TOPSIS method

In the TOPSIS approach, the alternative with shortest negative distance and the farthest positive distance from the ideal solution is preferred. The TOPSIS solution method is performed in a five stage algorithm as follows:

a. Normalize the decision matrix as follows:

$$n_{ij} = \frac{r_{ij}}{\sqrt{\sum_{j=1}^{n} r_{ij}^2}} \quad i=1,2,...,m \quad j=1,2,...,n \tag{2.7}$$

b. Assume that $w_i \quad i=1,2,...,m$ are the weight of the *i*th attribute. Then,

$$V_{ij} = n_{ij} w_i \quad i=1,2,...,m \quad j=1,2,...,n \tag{2.8}$$

c. Assume that $K$ and $K'$ are the index sets of cost and benefit attributes, respectively. Then, Find the positive and negative ideal solutions as follows

$$\text{Positive ideal solution } V^+ = \{(\max_j V_{ij} \mid i \in K'), (\min_j V_{ij} \mid i \in K)\} \tag{2.9}$$
$$= \{V_1^+, V_2^+, ..., V_l^+\}$$
$$\text{Negative ideal solution } V^- = \{(\min_j V_{ij} \mid i \in K'), (\max_j V_{ij} \mid i \in K)\}$$
$$= \{V_1^-, V_2^-, ..., V_l^-\}$$

d. Measure the distances from ideal solution.

$$S_j^+ = \left( \sum_{i=1}^{m} (V_{ij} - V_i^+)^2 \right)^{0.5} \quad j=1,2,...,n \tag{2.10}$$

$$S_j^- = \left( \sum_{i=1}^{m} (V_{ij} - V_i^-)^2 \right)^{0.5} \quad j=1,2,...,n$$

e. Calculate the relative closeness to the ideal solution.

$$C_j = \frac{S_j^-}{S_j^+ + S_j^-} \tag{2.11}$$

## 3. New approach

### 3.1 Criticisms the ranking of DEA models

Purposes of DEA rely on two issues, ranking and recommendation for improving inefficient DMUs by decreasing inputs and/or increasing outputs. Several DEA papers have developed ranking methods, and attempted to discriminate DMUs in order to determine priority of them. In spite of attempting for developing ranking method, there is no study to investigate the usage of the results derived from ranking. The examples below illustrate how unreliable the ranking methods are, and that using DEA models can mislead researchers.

#### 3.1.1. Comparison the priority of an inefficient DMU and an efficient one

To criticize the validity of ranking in DEA, we refer to an example. At first, two efficient and inefficient DMUs are compared. The following table includes 4 DMUs with single input and single output.

**Table 3.1.** Data of four DMUs with single input and output for criticism of DEA ranking

|            | $DMU_1$ | $DMU_2$ | $DMU_3$ | $DMU_4$ |
|------------|---------|---------|---------|---------|
| Input      | 1       | 2       | 2       | 2       |
| Output     | 1       | 2       | 10      | 9       |
| Efficiency | 1       | 0.556   | 1       | 0.944   |

The figure below illustrates the locus of the above DMUs. Using the concept of return to scale, it is known that Increasing Return to Scale (IRS) is proper to evaluate the above DMUs.

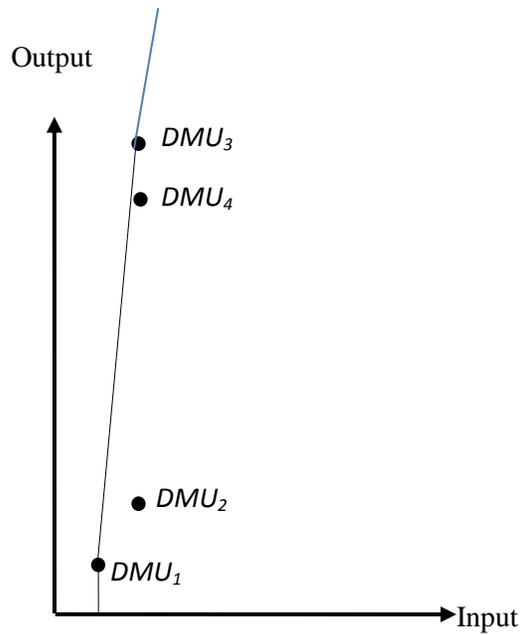

**Figure3.1.** DEA frontier corresponding to four DMUs

The efficiency scores using IRS model are seen in the last row of Table 1. As shown, $DMU_1$ and $DMU_3$ are considered as efficient and both DMUs 2 and 4 are inefficient. If it is required to choose two the best DMUs amongst them, which ones are preferred? To clarify the issue, a real world problem is matched on the example. Assume that DMUs in the above example are investment companies, and an investor would like to utilize his capital and invest his money in two companies. Input and output are assumed to be the amount of investment (in million dollar) and the profit after 5 years (in million Dollar), respectively. In other words, investor would like to rank the priority of his options for investment. Using DEA, DMUs 1, 3, 4 and 2 are placed in ranks 1 to 4, respectively, but logically investment in project 4 has more benefit than project 1. By applying Max-min and TOPSIS different results are obtained, and they confirm that the DMUs 4 and 3 are the best two options for investment.

### 3.1.2. Comparison the priority of efficient DMUs

Furthermore, let us analyze the most well-known super efficiency model, say AP. To evaluate the validity of this model, we refer to the above example. First, we eliminate $DMU_4$ and evaluate two efficient DMUs to rank them. Efficiency scores of DMUs 1 and 3 obtained from Model (2) by omitting the DMU under evaluation are 2 and 5, respectively. These scores confirm that $DMU_3$ is preferred to $DMU_1$ which is reliable. Again, we include $DMU_4$ and evaluate efficient DMUs 1 and 3 using AP model. The efficiency scores of DMUs 1 and 3 will be 2 and 1.11, and consequently, order of priority of DMUs is changed. In fact, $DMU_1$ is recommended as the best option. What changes the order of DMUS is the existence of a powerful inefficient DMU$_4$ which places close to DMU$_3$. This closeness influences on the efficiency of DMU$_3$. Undoubtedly, DMUs 3 and 4 have got better performance than DMU$_1$. As it is seen, existence of DMU$_4$ worsens the status of the efficient DMU$_3$. This is marvelous that the existence of an inefficient DMU can change the priority of DMUs. This means that choosing project one is the best option for investment, but logically it cannot be preferred to DMU$_3$. Reference set analysis attempts to remedy this shortcoming.

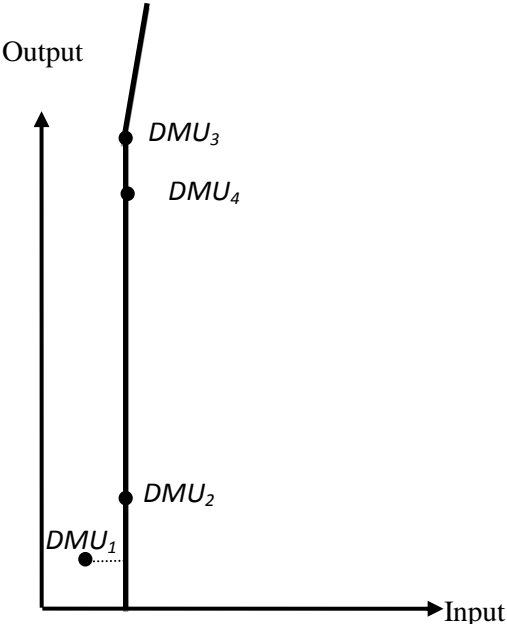 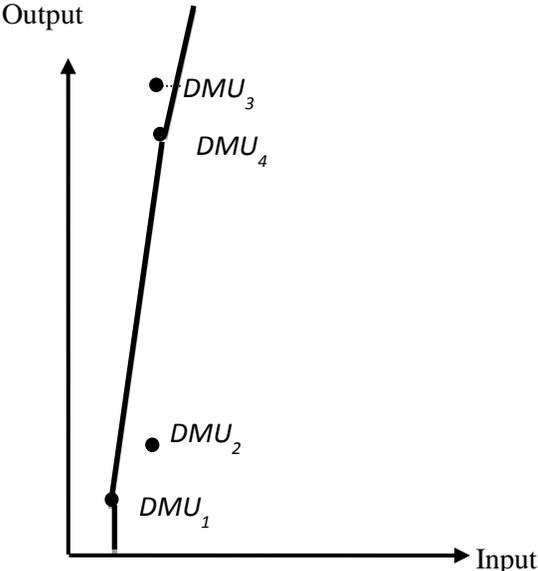

**Figure3.1.a**  **Figure3.1.b**

DEA frontiers corresponding to four DMUs

To evaluate the performance of reference set analysis approach, few more DMUs are added to PPS. Given that DMUs 5,6,7,8, and 9 are included to the evaluation, and their positions are seen in the following figure:

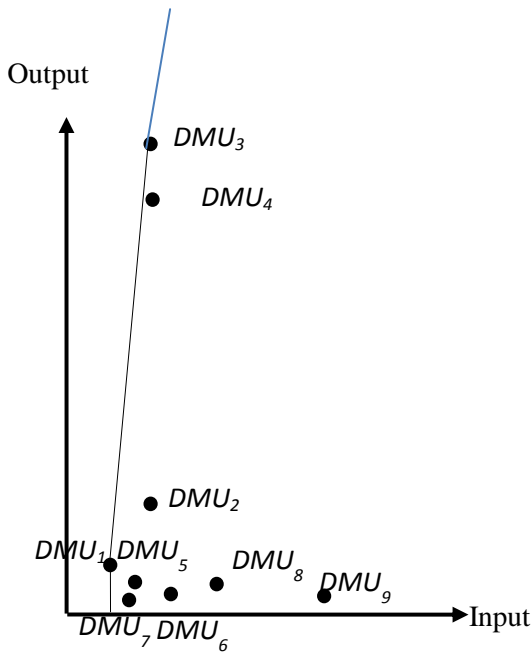

**Figure3.2.** DEA frontier corresponding to nine DMUs

In view of reference set analysis, being reference set of 5 inefficient DMUs, $DMU_1$ has higher rank than $DMU_3$. Based on this view, two powerful DMUs 3 and 4 are apart from others because of their super-efficient and outstanding performance. Performing in higher level than others, DMUs 3 and 4 are being criticized, and other DMUs actually are not able to perform in such a high level. On other words, inefficient and powerless DMUs 2, 5, 6, 7, 8, and 9 justify the performance of the weak $DMU_1$.

As a result, the indicated above examples substantiate that the DEA models are not reliable for ranking. Ultimately, the question that remains in this research is "If we assume that DEA ranking methods are not to determine priority of DMUs, what are they for?"

### 3.2 Hybrid self-assessment DEA model (SADEA) basis on fuzzy concept

The shortcoming presented in subsection 3.1 implies that DEA models are unable to rank DMUs efficiently, and a powerful model is required to do that. In this section, a novel methodology basis on fuzzy concept is presented in which the DMU under evaluation is appraised as a self-assessment methodology. The method is called self-assessment because we do not use other DMUs' data, except the maximum value of output to create the fuzzy number corresponding to output. In fact, each DMU gains its efficiency score considering distance from the best DMU and some fuzzy numbers indicating the satisfaction of DM (decision maker). The mentioned satisfaction depends on various factors. For instance, increasing products from the current level is possible with force or encourage of workers through payment, so, this involves expenses. On the other hand, the satisfaction of decision maker is directly related to difficulties in increasing production. The more difficulties in increasing production the less satisfaction will be caused for DM. For a better understanding of the models which will present in future, we refer to Figure (3.1) which provides a two-dimensional diagram of a simple efficiency case study involving 5 DMUs (A, B, C, D and E) in which only a single input is used to produce a single output.

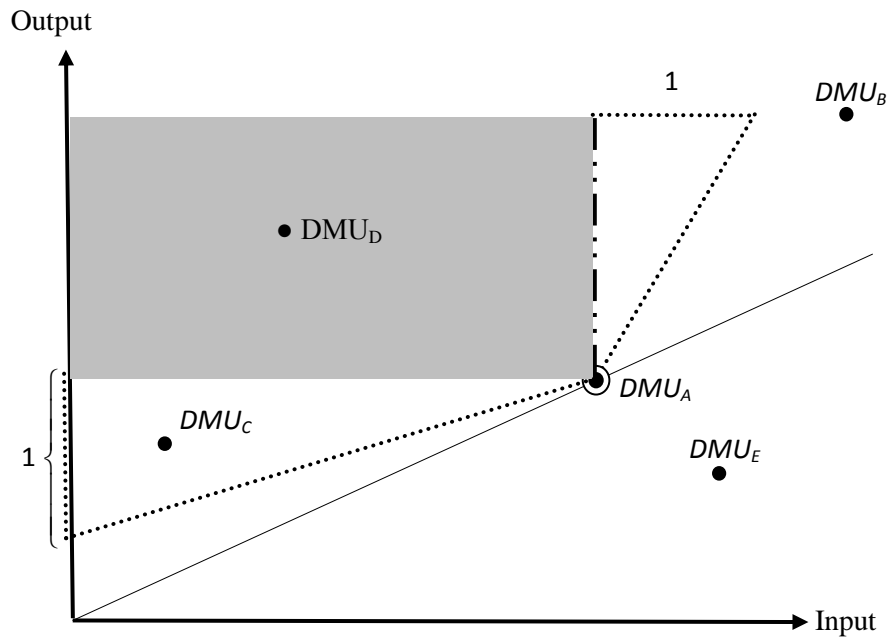

**Figure3.2.** Effective Space highlighted in this figure

The shaded space in figure 3.2 demonstrates a space in which $DMU_A$ is able to improve its efficiency by decrease in the input and increase in the output. Hereinafter, the above mentioned surface is called Effective Space. The DMU under evaluation categorizes other DMUs in two parts; 1) DMUs that affect its inefficiency improvement 2) DMUs that do not affect its inefficiency improvement. The efficiency frontier of each DMU determines the frontier between two parts. For each DMU in the part 1, there is a DMU in the effective space with the same efficiency, because efficiency frontier of each DMU in part 1 passes through effective space. The connecting line between the points corresponding to DMUs C and D passes through origin of coordinates. Therefore, $DMU_D$ is the representative of $DMU_C$ in the effective space. On the other words, each point on the line passing through $DMU_C$ and origin of coordinates located in the effective space is a representative of $DMU_C$. Dash lines in Figure 3.2 illustrate the membership function of the fuzzy numbers related to input and output. These membership functions are super imposed onto the graph.

In the proposed model we are looking for a point (a DMU) in the shaded space with further distance from the efficiency frontier of the DMU under evaluation. So, the evaluation in this methodology is kind of a super efficiency evaluation, however there is not any specific DMU. On the other words, each point in the shaded space is playing the role of a DMU under evaluation. As it is seen in figure 3.3, there is relationship between the slopes of the lines associated with the membership function and closeness to the coordinate's axes. The input value is much larger, gentler slope of the line will be.

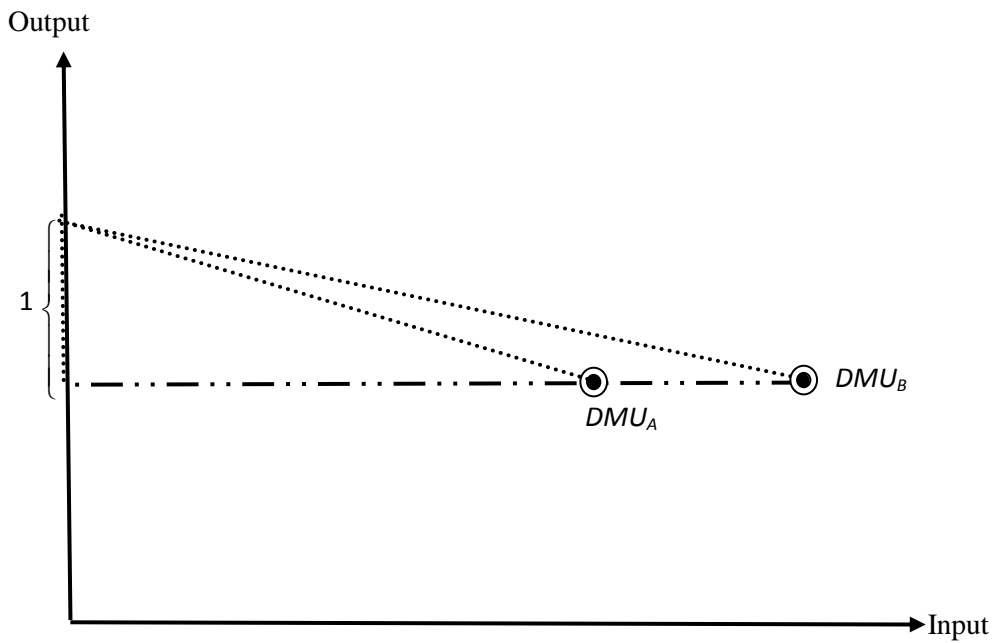

**Figure3.3.** The more slopes of the membership function the closer to the coordinate's axes

Based on the mentioned above explanation, the following fuzzy linear programming is proposed evaluating $DMU_p$:

$$\min \quad w_p$$
$$\text{s.t.} \tag{3.1}$$
$$\lambda x_{ip} \leq \tilde{x}_{ip} + w_p \quad \forall i$$
$$\lambda y_{rp} \geq \tilde{y}_{rp} - w_p \quad \forall r$$
$$\lambda \geq 0$$
$$w_p \text{ free}$$

Model (3.1) is a fuzzy linear programming with fuzzy resources. As it is seen, there are just two DMUs in Model (3.1), the DMU under evaluation $p$ and a fuzzy DMU with the input vector $\tilde{X}_p = (\tilde{x}_{1p}, \tilde{x}_{2p}, ..., \tilde{x}_{mp})^t$ and the output vector $\tilde{Y}_p = (\tilde{y}_{1p}, \tilde{y}_{2p}, ..., \tilde{y}_{sp})^t$. Actually, $DMU_p$ is being appraised independently. In fact, in Model (3.1) the fuzzy DMU with the input and output vectors $\tilde{X}_p$ and $\tilde{Y}_p$ is being evaluated. Since our purpose is to evaluate $DMU_p$, therefore, the efficiency score is larger; the $DMU_p$ is placed in the lower rank.

In Model (3.1), the membership values of the fuzzy numbers $\tilde{x}_{ip}$ and $\tilde{y}_{rp}$ are defined as follows:

$$\mu_{\tilde{x}_{ip}}(\bar{x}_{ip}) = \frac{\bar{x}_{ip} - 0}{x_{ip} - 0} \quad \text{for } 0 \leq \bar{x}_{ip} \leq x_{ip} \quad (i = 1, 2, ..., m) \tag{3.2}$$

$$\mu_{\tilde{y}_{rp}}(\bar{y}_{rp}) = \frac{\bar{y}_{rp} - y_{rp}}{y_r^{\max} - y_{rp}} \quad \text{for } y_{rp} \leq \bar{y}_{rp} \leq y_r^{\max} \quad (r = 1, 2, ..., s) \tag{3.3}$$

in which $y_r^{\max} = \max_j \{y_{rj}\}$.

Since the inputs and outputs are not homogeneous and scale of objective function in the proposed model (3.1) depends on the units of measurement of inputs and outputs data, unit dependence is obtained by normalization e.g. dividing each input to the largest of them and each output to the largest of them as one of the techniques for normalization. For fuzzy data consider the following definition:

**Definition**: consider the fuzzy numbers $\tilde{x}_i = (x_i^m, x_i^u)$ $(i = 1, 2, ...k)$. We introduce $\tilde{x}_i' = (x_i'^m, x_i'^u)$ $(i = 1, 2, ...k)$ as follows and we call them normalized $\tilde{x}_i$:

$$\mu_{\tilde{x}'_i}(\overline{x}'_i) = \mu_{\tilde{x}_i}(\overline{x}_i) \text{ and } \overline{x}'_i = \frac{\overline{x}_i}{\max_i\{x_i^u\}} \quad (x_i^m \leq \overline{x}_i \leq x_i^u) \quad (i = 1, 2, \ldots, k)$$

Assume that the fuzzy numbers $\tilde{x}'_{ip}$ $(i = 1, 2, \ldots m)$ and interval outputs $\tilde{y}'_{rp}$ $(r = 1, 2, \ldots s)$ are the normalized numbers of $\tilde{x}_i$ and $\hat{y}_r$. In addition, assume $x'_{ip} = \frac{x_{ip}}{\max_i\{x_{ip}\}}$ and $y'_{rp} = \frac{y_{rp}}{\max_r\{r_{rp}\}}$.

Model (3.1) is a super efficiency model. To solve the fuzzy linear programming (3.1), it is suggested the following multi objective linear programming:

$$\max \quad \mu_{\tilde{x}'_{ip}}(\overline{x}'_{ip}) \tag{3.4}$$
$$\max \quad \mu_{\tilde{y}'_{rp}}(\overline{y}'_{rp})$$
$$\min \quad w_p$$
s.t.
$$\lambda x'_{ip} \leq \overline{x}'_{ip} + w_p \quad \forall i$$
$$\lambda y'_{rp} \geq \overline{y}'_{rp} - w_p \quad \forall r$$
$$0 \leq \overline{x}'_{ip} \leq x'_{ip}$$
$$y'_{rp} \leq \overline{y}'_{rp} \leq \frac{y_r^{\max}}{y_r^{\max}} = 1$$
$$\lambda \geq 0$$
$$w_p \text{ free}$$

Model 3.4 is converted to linear programming problem as follows:

$$\max \quad \alpha \qquad (3.5)$$

s.t.

$$\lambda x'_{ip} \leq \overline{x}'_{ip} + w_p \quad \forall i$$
$$\lambda y'_{rp} \geq \overline{y}'_{rp} - w_p \quad \forall r$$
$$0 \leq \overline{x}'_{ip} \leq x'_{ip} \quad \forall i$$
$$y'_{rp} \leq \overline{y}'_{rp} \leq 1 \quad \forall r$$
$$\alpha \leq \mu_{\tilde{x}'_{ip}}(\overline{x}'_{ip}) \quad \forall i$$
$$\alpha \leq \mu_{\tilde{y}'_{rp}}(\overline{y}'_{rp}) \quad \forall r$$
$$\alpha \leq 1 - w_p$$
$$\lambda \geq 0$$
$$w_p \text{ free}$$

As mentioned above, there is an expressive relationship between the distance from inputs and outputs to coordinates exes and the slope of lines corresponding to membership functions. The steep slope causes more membership value with input decreasing or output increasing. The result of the new ranking method proposed in this paper is essentially different from traditional DEA models. For example, despite the traditional DEA model which introduces projects 1 and 3 as two the best options for investment, by running model (3.5), projects 3 and 4 are recommended in Example1. The following example illustrates the logic of the new approach.

**Example2.** Consider 4 DMUs with single input and single output listed in Table 3.2.

**Table3.2.** Data of four DMUs with single input and output

|  | $DMU_1$ | $DMU_2$ | $DMU_3$ | $DMU_4$ |
|---|---|---|---|---|
| Inputs | 1 | 1 | 2 | 6 |
| Outputs | 8 | 7 | 2 | 1 |

The effective spaces corresponding to the above DMUs have been illustrated in the Figures 3.4 to 3.7. As it is seen, $DMU_4$ and $DMU_1$ have the most and the least spaces, respectively, compared with two other DMUs. The slopes of the membership function related to $DMU_4$ are gentler than that of other DMUs. Therefore, in this evaluation, we expect that $DMU_4$ gains the least efficiency score. With the same reduction in the amount of inputs, $DMU_1$ and $DMU_2$ achieve the larger membership values compared with $DMU_3$ and $DMU_4$. So, there is a relationship between this superiority and inputs being small.

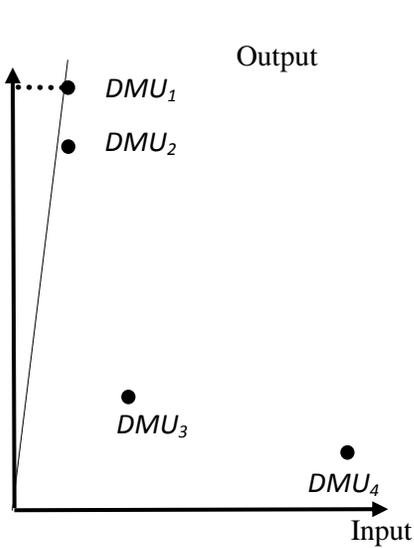

**Figure 3.4.** Effective space of $DMU_1$

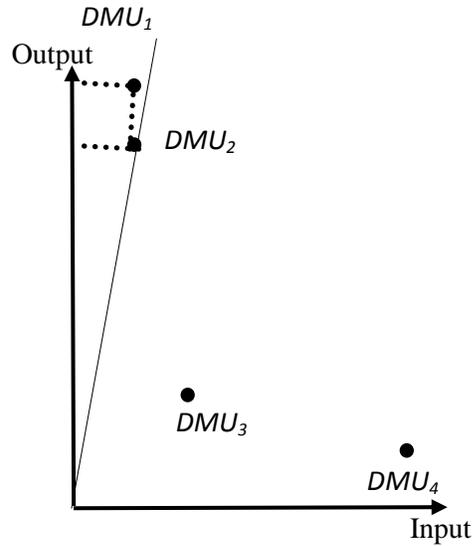

**Figure 3.5.** Effective space of $DMU_2$

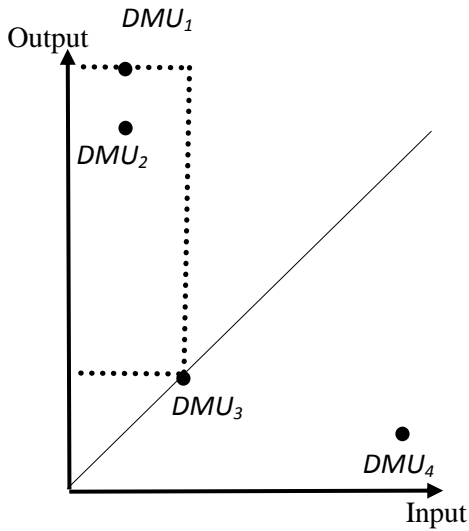

**Figure 3.6.** Effective space of $DMU_3$

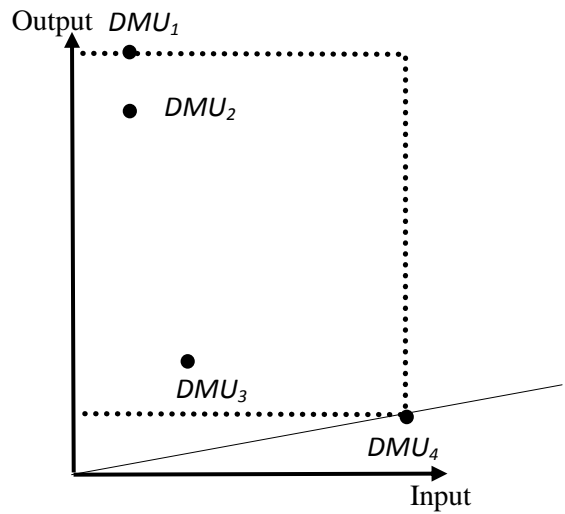

**Figure 3.7.** Effective space of $DMU_4$

The efficiency scores obtained by Model (3.5) are demonstrated in Table (3.3)

**Table3.3.** Evaluation of 4 DMUs with single input and output using proposed method

|  | $DMU_1$ | $DMU_2$ | $DMU_3$ | $DMU_4$ |
|---|---|---|---|---|
| Efficiency | 0.1250000 | 0.1379310 | 0.3636364 | 0.4705882 |

In the proposed method, the efficiency score is smaller; the higher evaluated DMU will be ranked. Therefore DMUs 1, 2, 3 and 4 are placed in ranks 1 to 4, respectively.

As it is seen, the result of efficiency scores in the CCR and proposed model is same, but it must not mislead us. Actually, the traditional DEA models and the proposed approach are different in nature. Although the new approach has been designed basis on a DEA model, there is more affinity among the proposed method and the MADM methodologies. One of the most important features of this methodology is that the model evaluates DMUs (or Alternatives) in a self-assessment system. As only the endpoints of the data interval play the role in evaluating alternatives in TOPSIS and Max-min, in the new methodology just the maximum outputs play such a role. In the proposed method a relative ideal is sought in the effective space. Despite TOPSIS, it is not required to weights in the proposed model, though it can be considered in. Furthermore, Max-min does not need information from decision makers. So, we will compare the proposed method with two latter i.e. the Max-min and TOPSIS methods. In TOPSIS the weights are considered as equal.

**4. Case Study**

As mentioned, Max-min does not need information from decision makers. Since the proposed methodology has such a similar structure, it is compared with these two methods. For this purpose, we implement the new methodology, Max-min and TOPSIS on the case study presented by Sowlati and Paradi [12]. Among 79 DMUs and 6 inputs and outputs, eighteen DMUs (alternatives) and five Inputs and outputs were selected. Inputs and outputs are costs and benefits, respectively. Two types of full time equivalent number of employees (FTE sales and FTE supports) were considered the inputs (cost alternatives) of the model. Loans, mortgages, registered retirement saving plans (RRSPs) and letters of credit (LC) were considered as the outputs (benefit alternatives) of the model. Table 3.4 demonstrates the data selected.

**Table 3.4.** Data collected by Sowlati and Paradi [12]

| Alternative | FTE sales | FTE support | RRSP | LC | Loans | Mortgages |
|---|---|---|---|---|---|---|
| 1 | 45.34 | 40.93 | 263 | 137 | 935 | 429 |
| 2 | 9.02 | 1.34 | 42 | 6 | 176 | 32 |
| 3 | 26.12 | 8.24 | 130 | 20 | 679 | 101 |
| 4 | 10.94 | 4.87 | 134 | 37 | 437 | 80 |
| 5 | 49.52 | 32.28 | 308 | 46 | 726 | 227 |
| 6 | 10.82 | 1.09 | 27 | 2 | 18 | 136 |
| 7 | 11.52 | 1.98 | 44 | 5 | 337 | 47 |
| 8 | 8.11 | 3.91 | 34 | 1 | 245 | 33 |
| 9 | 9.96 | 5.26 | 29 | 2 | 202 | 49 |
| 10 | 9.86 | 1.01 | 67 | 10 | 161 | 52 |
| 11 | 4 | 1.58 | 42 | 2 | 159 | 17 |
| 12 | 5.78 | 1.52 | 85 | 1 | 196 | 78 |
| 13 | 4.87 | 1.05 | 52 | 4 | 237 | 52 |
| 14 | 2.93 | 1.97 | 6 | 2 | 127 | 18 |
| 15 | 2.96 | 1.58 | 21 | 2 | 103 | 23 |
| 16 | 9.84 | 5.02 | 55 | 1 | 301 | 50 |
| 17 | 16.06 | 1.99 | 143 | 7 | 551 | 187 |
| 18 | 25.06 | 7.76 | 151 | 13 | 808 | 211 |

The result of the Max-min ranking method is seen in Table 3.5. In Max-min, priority of alternatives located in ranks 1 to 4 is affected by the number 137 related to alternative 1. The mentioned number is very large compared with other numbers associated to the output LC. This causes that, after normalization, the minimum of the alternatives 3, 4 and 18 are selected. Therefore, the Min-max ranking method is not reliable.

**Table 3.5.** The result of Max-min method

|  | Scores 1st level | Scores 2nd level | Scores 3rd level | Scores 4th level | Max-min Rank |
|---|---|---|---|---|---|
| 1 |  |  | 0.025654 |  | 10 |
| 2 |  | 0.043796 |  |  | 6 |
| 3 | 0.112175 |  |  |  | 2 |
| 4 | 0.186480 |  |  |  | 1 |
| 5 |  | 0.032528 |  |  | 8 |
| 6 |  |  |  | 0.022277 | 15 |
| 7 |  | 0.036496 |  |  | 7 |
| 8 |  |  |  | 0.021739 | 16-17-18 |
| 9 |  |  |  | 0.043478 | 11-12-13 |
| 10 |  | 0.072993 |  |  | 4 |
| 11 |  |  |  | 0.043478 | 11-12-13 |
| 12 |  |  |  | 0.021739 | 16-17-18 |
| 13 |  |  | 0.029197 |  | 9 |
| 14 |  |  |  | 0.022814 | 14 |
| 15 |  |  |  | 0.043478 | 11-12-13 |
| 16 |  |  |  | 0.021739 | 16-17-18 |
| 17 |  | 0.051095 |  |  | 5 |
| 18 | 0.094891 |  |  |  | 3 |

The results of the proposed method and TOPSIS are seen in Table 3.6. In this case, the weights of the alternative have been considered as equal.

**Table 3.6.** The results of the proposed and TOPSIS

| Alternative | Scores in the new method | Rank | TOPSIS score | Rank |
|---|---|---|---|---|
| 1 | 0.3503982 | 1 | 0.597509 | 1 |
| 2 | 0.4463411 | 10 | 0.374171 | 14 |
| 3 | 0.4392190 | 6 | 0.398698 | 6 |
| 4 | 0.3516046 | 2 | 0.450643 | 2 |
| 5 | 0.4281250 | 4 | 0.439257 | 5 |
| 6 | 0.4838366 | 15 | 0.379169 | 10 |
| 7 | 0.4706567 | 13 | 0.376565 | 12 |
| 8 | 0.4891005 | 17 | 0.363366 | 17 |
| 9 | 0.4824899 | 14 | 0.352735 | 18 |
| 10 | 0.4225485 | 3 | 0.385865 | 9 |
| 11 | 0.4585621 | 11 | 0.378927 | 11 |
| 12 | 0.4848400 | 16 | 0.396415 | 7 |
| 13 | 0.4353719 | 5 | 0.392661 | 8 |
| 14 | 0.4450913 | 7 | 0.372725 | 15 |
| 15 | 0.4455872 | 8 | 0.376274 | 13 |
| 16 | 0.4909822 | 18 | 0.364603 | 16 |
| 17 | 0.4634891 | 12 | 0.448561 | 4 |
| 18 | 0.4571410 | 9 | 0.449703 | 3 |

## 5. Conclusion

In this paper, each evaluation includes two DMUs; the DMU under evaluation and a fuzzy DMU which is representative of the effective space. In this space DMUs can maneuver to improve their situation. To this purpose, by considering the difficulty of increasing outputs by retaining the level of current inputs and also the difficulty of decreasing inputs by retaining the level of current outputs, some fuzzy numbers were introduced. The mentioned fuzzy numbers together formed a fuzzy DMU. Each DMU was evaluated by its own fuzzy DMU different from others. In fact, the efficiency of DMUs were obtained basis on the potential of inputs and outputs to decrease and increase, respectively. This type of assessment would be classified as a multi-attribute decision making method. Furthermore, despite TOPSIS, which applies just maximum and minimum data in the process of performance evaluation, the proposed model, by introducing fuzzy numbers, allows all data to cooperate in evaluation. As discussed, the effective space is narrower; the DMU under evaluation is more efficient. The contribution of this method is that each DMU is evaluated as self-assessment. In each evaluation, DMUs hunt for their ideal DMU, independently.

## References

[1] A. Charnes, W.W. Cooper, E. Rhodes, Measuring the efficiency of decision making units, *Eur. J. Oper. Res.* 2 (1978) 429–444.

[2] R.D. Banker, A. Charnes, W.W. Cooper, Some methods for estimating technical and scale inefficiencies in data envelopment analysis, *Manage. Sci.* 30 (1984) 1078–1092.

[3] N. Adler, L. Friedman, Z. Sinuany-Stern, Review of ranking methods in the data envelopment analysis context, *Eur. J. Oper. Res.* 140 (2002) 249–265.

[4] A. Charnes, W.W. Cooper, B. Golany, L. Seiford, J. Stutz, Foundations of data envelopment analysis for Pareto–Koopmans efficient empirical production functions. *J. Econ.* 30 (1985) 91–107.

[5] K. Tone, A slack-Based measure of efficiency in data envelopment Analysis, *Eur. J. Oper. Res.* 130 (2001) 498-509.

[6] T. R. Sexton, R.H. Silkman, A. Hogan, Data envelopment analysis: Critique and extensions, In Silkman, R. H. (Ed), Measuring Eficiency: An Assessment of Data Envelopment Analysis,

Publication no. 32 in the series New Direction of Program Evaluation. San Francisco: Jossey Bass. (1986) 73–105.

[7] P. Andersen, N.C, Petersen, A procedure for ranking efficient units in data envelopment analysis. *Manage. Sci.* 39(1993) 1261–1264.

[8] S. Mehrabian, M. Alirezaee, G. Jahnshahloo, A Com-plete Efficiency Ranking of Decision Making Units in Data Envelopment Analysis, *Comput. Optim. Appl.* 14 (1999) 261-266.

[9] S. Saati, M. Zarafat Angiz, A. Memariani, G.R. Jahanshahloo, A model for ranking decision making units in data envelopment analysis, *Ricerca Operativa* 31(97) (2001) 47- 59.

[10] Z. Sinuany-Stern, A. Mehrez, Y. Hadad, An AHP/DEA methodology for ranking decision making, *Int. T. Oper. Res. 7* (2000) 109–124.

[11] D. Dubois, H. Prade, The mean value of a fuzzy number, *Fuzzy Set. Syst.* 24 (1987) 279–300.

[12] T. Sowlati, J.C. Paradi, Establishing the "practical frontier" in data envelopment analysis, *Omega.* 32 (2004) 261 – 272.